\documentclass[reqno,12pt]{amsart}

\usepackage{amsmath}
\usepackage{amscd}
\usepackage{amssymb}
\usepackage{enumerate}
\usepackage{amsfonts}

\DeclareFontEncoding{OT2}{}{}

\numberwithin{equation}{section}

\newcommand{\rar}[1]{\stackrel{#1}{\longrightarrow}}

\newcommand{\Ga}{\Gamma}  
\newcommand{\la}{\lambda}

\newcommand{\bC}{{\mathbb C}}

\newcommand{\bF}{{\mathbb F}}

\newcommand{\bN}{{\mathbb N}}

\newcommand{\bZ}{{\mathbb Z}}

\newcommand{\Gh}{\widehat{G}}

\newcommand{\Kh}{\widehat{K}}

\newcommand{\gh}{\widehat{\Ga}}

\newcommand{\abs}[1]{\vert #1\vert}

\newcommand{\id}{\operatorname{id}}

\newcommand{\Ind}{\operatorname{Ind}}
\newcommand{\Res}{\operatorname{Res}}

\newcommand{\tens}{\otimes}

\newcommand{\sbr}{\smallbreak}
\newcommand{\mbr}{\medbreak}

\newtheorem*{thmnn}{Theorem}

\newtheorem{thm}{Theorem}[section]
\newtheorem{cor}[thm]{Corollary}
\newtheorem{lem}[thm]{Lemma}
\newtheorem{prop}[thm]{Proposition}

\theoremstyle{remark}

\newtheorem{defin}[thm]{Definition}

\newcommand{\cst}{\bC^\times}

\newcommand{\Fr}{\operatorname{Fr}}

\newcommand{\Gal}{\operatorname{Gal}}

\newcommand{\fdim}{\operatorname{fdim}}

\newcommand{\sh}{\operatorname{sh}}

\newcommand{\sheis}{\operatorname{SHeis}}

\newcommand{\ti}[1]{\widetilde{#1}}

\title[Base change maps for unipotent algebra groups]{Base change maps for unipotent algebra groups}

\author[Mitya Boyarchenko]{Mitya Boyarchenko}

\thanks{The author's research was partially supported by NSF grant
DMS-0401164.}

\begin{document}

\maketitle

\setcounter{tocdepth}{1}

\tableofcontents

\section*{Introduction}

Let $p$ be a fixed prime, and let $\bF$ be a fixed algebraic
closure of a finite field with $p$ elements. If $q$ is a power of
$p$, we will denote by $\bF_q$ the unique subfield of $\bF$ that
consists of $q$ elements. We will denote by $\Fr_q$, or simply
$\Fr$ if no confusion is possible, the canonical topological
generator of the Galois group $\Gal(\bF/\bF_q)$, given by
$x\mapsto x^q$. By a {\em representation} of a finite group we
will always mean a complex finite dimensional representation.

\mbr

Let $G$ be a unipotent algebraic group over $\bF_q$. For every
$n\in\bN$, let $\Ga_n=G(\bF_{q^n})$, a finite nilpotent group. We
have the obvious action of $\Gal(\bF_{q^n}/\bF_q)$ on $\Ga_n$, and
hence an induced action of $\Gal(\bF_{q^n}/\bF_q)$ on the set $\gh_n$ of
isomorphism classes of its irreducible representations.
The three basic questions in geometric character theory for
unipotent groups over finite fields are as follows (cf.
\cite{drinfeld}, which was motivated in part by \cite{lusztig}):

\mbr

\begin{enumerate}[(1)]
\item Do there exist ``base change maps''
\[
T_m^n : \gh_m \rar{} \gh_n^{\Gal(\bF_{q^n}/\bF_{q^m})}
\]
for all $m,n\in\bN$ such that $m\lvert n$, which are
$\Fr$-equivariant and satisfy
\[
T_m^k = T_n^k \circ T_m^n
\]
for all $m,n,k\in\bN$ such that $m\lvert n\lvert k$\,? \sbr
\item Assuming that the answer to question (1) is positive,
are the base change maps injective (resp., surjective; resp.,
bijective)? \sbr
\item Assuming again that the answer to question (1) is positive,
form the direct limit
\[
\gh := \underset{\rar{}}{\lim} \gh_n
\]
with respect to the base change maps, and equip it with the
induced action of $\Gal(\bF/\bF_q)$. Does there exist a
``geometric object'' $\Gh$ defined over $\bF_q$ such that its set
of geometric points $\Gh(\bF)$ admits an $\Fr$-equivariant
bijection onto $\gh$?
\end{enumerate}

\sbr

\noindent
 We prefer not to make the last question precise at this point.

\mbr

In this paper we study the first two questions for a special class
of groups. We recall that a {\em finite algebra group} (cf.
\cite{isaacs}) is a group of the form $1+A$, where $A$ is a finite
dimensional associative nilpotent algebra over $\bF_q$, and $1+A$
consists of formal expressions of the form $1+a$, where $a\in A$,
multiplied in the natural way: $(1+a)(1+b)=1+(a+b+ab)$. It is
clear that one can also think of $G=1+A$ as an affine algebraic
group over $\bF_q$: the underlying variety is the affine space
over $\bF_q$ underlying the algebra $A$, and for any commutative
$k$-algebra $R$, we have $G(R)=1+R\tens_k A$. When thought of in
this way, we will call $G$ a {\em unipotent algebra group} to
emphasize that we consider it as a geometric, rather than
algebraic, object. We will prove the following

\begin{thmnn}
There exist canonically defined injective base change maps for
irreducible representations of finite algebra groups.
\end{thmnn}

\mbr

The main technical tool used in the proof is the notion of a
``strongly Heisenberg representation'' of a finite algebra group,
introduced below. Apart from its significance in the geometric
theory, this notion also seems to have applications to purely
algebraic questions: for example, it can be used to rewrite
Halasi's proof \cite{halasi} of Gutkin's conjecture \cite{gutkin}
(see Section \ref{s:gutkin-halasi}) in such a way that it no
longer uses Isaacs' Theorem A \cite{isaacs}.

\mbr

It is interesting to note that the only previously known examples
of base change maps for unipotent groups over finite fields arise
from situations where there exists an explicit classification of
irreducible representations of the groups $\Ga_n$ in the spirit of
the orbit method and/or Mackey's little groups method. For
instance, an analogue of Kirillov's orbit method bijection can be
used to define natural base change maps when $G$ is a connected
unipotent group over $\bF_q$ that has nilpotence class less than
$p$. In all such situations, one effectively starts from a
positive answer to question (3), and uses it to answer questions
(1) and (2). On the other hand, question (3) for unipotent algebra
groups is still open, and no explicit classification of irreducible representations is known. It is also worth noting that in the
situations where the orbit method does work, the base change maps
one obtains are necessarily surjective, but may fail to be
injective. On the other hand, the base change maps constructed in
our paper are always injective. We do not know if they are always
surjective, but we expect that the answer is negative, even though
it is positive in many special cases.

\subsection*{Acknowledgements} I am deeply grateful to Vladimir
Drinfeld for introducing me to the basic problems of geometric
character theory mentioned above; for explaining to me a canonical
way of writing any irreducible representation of a finite
nilpotent group as the induced representation of a Heisenberg
representation, which has led to the reduction process described
in Section \ref{s:reduction}; and for the suggestion that the norm
maps introduced in Section \ref{s:base-1-dim} may be defined using
the functor $K_1$. I would also like to thank Charles Weibel for
answering my questions about $K$-theory and, in particular, for
his help in proving part (4) of Theorem \ref{t:norms}.

\section{Strongly Heisenberg
representations}\label{s:strongly-heisenberg}

\subsection{Heisenberg representations}
Once and for all, we make the convention that round parentheses
$(\cdot,\cdot)$ denote the group commutator, while square brackets
$[\cdot,\cdot]$ denote the Lie bracket (which in particular is
well defined in any associative algebra). In order to motivate the
definition of a strongly Heisenberg representation, we first
recall the more familiar notion of a Heisenberg representation. We
say that an irreducible representation $\rho:G\to GL(V)$ of a
finite group $G$ is a {\em Heisenberg representation} if there
exists a subgroup $H\subseteq G$ such that $G/H$ is abelian and
$H$ acts on $V$ by scalars, i.e., $\rho(H)\subseteq \cst$. Of
course, if $\rho$ is Heisenberg, then it is clear that one can
take $H=(G,G)$ in the definition, but it is often convenient to
allow $H$ to be larger. It is easy to check that the following conditions
are equivalent:
\begin{enumerate}[(i)]
\item $\rho$ is Heisenberg;
\item the image $\rho(G)\subseteq GL(V)$ is a nilpotent group of nilpotence class at
most $2$;
\item $G/Z(\rho)$ is abelian, where the center $Z(\rho)$ of $\rho$
is defined as
\[
Z(\rho) = \{ g\in G \big\lvert (\rho(g),\rho(h))=1 \ \forall\,
h\in G\} = \{ g\in G \big\lvert \rho(g)\in\cst \}.
\]
\end{enumerate}

\mbr

We now recall a basic fact. Let $G$ be an arbitrary finite group,
let $\rho : G \rar{} GL(V)$ be a Heisenberg representation, and
let $H\subseteq G$ be a subgroup such that $G/H$ is abelian and
$\rho(H)\subseteq\bC^\times$. The following result is well known
(and easy to prove).

\begin{lem}\label{l:heisenberg}
The commutator pairing
\[
c_\rho : G/H \times G/H \rar{} \bC^\times
\]
given by
\[
c_\rho(\overline{g}_1,\overline{g}_2)\cdot \id_V = \rho(g_1 g_2
g_1^{-1} g_2^{-1})
\]
is well defined and yields an alternating bilinear form on the
finite abelian group $G/H$. If $L\subseteq G/H$ is any maximal
isotropic subgroup with respect to $c_\rho$, and
$\widetilde{L}\subseteq G$ is its preimage in $G$, then $\rho$ is
induced from a one-dimensional representation of $\widetilde{L}$.
Furthermore, $\rho$ is determined, up to isomorphism, by the
character via which it acts on $Z(\rho)$.
\end{lem}

\subsection{Strongly Heisenberg representations} Let $k\subset\bF$
be a finite subfield, and let $G=1+A$ be a finite algebra group
over $k$.
\begin{defin}\label{d:strongly-Heisenberg}
An irreducible representation
\[
\rho : 1+A \rar{} GL(V)
\]
is said to be {\em strongly Heisenberg} if the subgroup
$1+A^2\subset 1+A$ acts by scalars on $V$, i.e.,
\[
\rho(1+A^2) \subseteq \bC^\times \subseteq GL(V).
\]
\end{defin}
Let $\rho$ be a strongly Heisenberg representation of a finite
algebra group $G=1+A$, and take $H$=$1+A^2$. It is clear that
$G/H$ is isomorphic to the underlying additive group of $A/A^2$,
and in particular is abelian, whence every strongly Heisenberg
representation is, a fortiori, Heisenberg. On the other hand,
$A/A^2$ is also a vector space over $k$, and we have the following
main result of the section:

\begin{prop}\label{p:strongly-Heisenberg}
In this situation, with the notation of Lemma \ref{l:heisenberg},
there exists a maximal isotropic subgroup $L\subseteq G/H=A/A^2$
with respect to $c_\rho$ which is also a $k$-subspace of $A/A^2$.
Consequently, $\ti{L}$ is an algebra subgroup of $G$.
\end{prop}

\mbr

The proof can be easily reduced to the following
\begin{lem}\label{l:strongly-Heisenberg}
Let $\phi:1+A^2\to\cst$ be a $G$-invariant character. The
commutator pairing
\[
c=c_\phi : (1+A)/(1+A^2) \times (1+A)/(1+A^2)  \rar{} \cst,
\]
given by $c(\overline{g},\overline{h})=\phi((g,h))$, is well
defined and satisfies the equation
\begin{equation}\label{e:commutator}
c(1+\la a,1+b) = c(1+a,1+\la b) \qquad \forall\, a,b\in A,\ \la\in
k.
\end{equation}
\end{lem}

Indeed, assume that this lemma has been proved. In the notation of
Proposition \ref{p:strongly-Heisenberg}, if $\phi$ is the
character by which $\rho$ acts on $1+A^2$, it is clear that
$c_\rho=c_\phi$. Now let $L\subseteq A/A^2$ be a $k$-subspace
which is maximal among all $k$-subspaces that are isotropic with
respect to $c_\rho$. Equation \eqref{e:commutator} implies that
the annihilator $L^\perp$ of $L$ with respect to $c_\rho$ is also
a $k$-subspace of $A/A^2$. If $L\neq L^\perp$, then
\eqref{e:commutator} also shows that $L$ together with any element
of $L^\perp\setminus L$ spans a $k$-subspace of $A/A^2$ which is
isotropic with respect to $c_\rho$, contradicting our choice of
$L$. This completes the proof of Proposition
\ref{p:strongly-Heisenberg}.

\subsection{Commutators in algebra groups}
In our proof of Lemma \ref{l:strongly-Heisenberg} we will need the
following result.
\begin{prop}[Halasi]\label{p:commutators} Let $R$ be an integral domain whose quotient
field has characteristic zero, and let $J=F_R(n,X)$ be a free
associative nilpotent algebra of class $n\in\bN$ generated by a
set $X$ over $R$. For all $k\geq 2$, we have
\[
(1+J,1+J)\cap (1+J^k) = (1+J,1+J^{k-1}).
\]
\end{prop}
Formally, $F_R(n,X)$ denotes the initial object in the category of
all maps from the set $X$ to associative nilpotent algebras $A$
over $R$ satisfying $A^n=(0)$. The proposition is proved in
\cite{halasi} in the special case $R=\bZ$; however, it is easy to
check that all steps in Halasi's proof work equally well for an
arbitrary integral domain of characteristic zero.

\subsection{Proof of Lemma \ref{l:strongly-Heisenberg}}
The $G$-invariance of $\phi$ is equivalent to
$\phi\bigl((1+A,1+A^2)\bigr)=\{1\}$. The fact that $c$ is well
defined is obvious from this. For the second statement of the
lemma, it suffices to prove that $\phi$ annihilates all products
of the form $(1+\la a,1+b)\cdot(1+a,1+\la b)^{-1}$. Now a
straightforward computation shows that $(1+\la
a,1+b)\cdot(1+a,1+\la b)^{-1}\in 1+A^3$, and, of course, we also
have $(1+\la a,1+b)\cdot(1+a,1+\la b)^{-1}\in (1+A,1+A)$. We would
be done if we could apply Proposition \ref{p:commutators} at this
point. In order to be able to use it, let $n\in\bN$ be such that
$A^n=(0)$, let $R$ be an arbitrary integral domain of
characteristic zero which admits a surjective homomorphism $R\to
k$, let $\widetilde{\la}\in R$ be a lift of $\la$, and let
$J=F_R(n,\{X,Y\})$ denote the free associative nilpotent algebra
of nilpotence class $n$ generated by two elements, $X$ and $Y$,
over $R$. Then, by Proposition \ref{p:commutators},
\[
(1+\widetilde{\la}X,1+Y)\cdot(1+X,1+\widetilde{\la}Y)^{-1}\in
(1+J,1+J)\cap (1+J^3) = (1+J,1+J^2).
\]
There exists a homomorphism of $R$-algebras $J\to A$ (where $A$ is
viewed as an $R$-algebra via $R\to k$) which takes $X\mapsto a$
and $Y\mapsto b$. Since it clearly takes $(1+J,1+J^2)$ into
$(1+A,1+A^2)$, the proof is complete.

\section{Theorems of Gutkin, Halasi, and Isaacs}\label{s:gutkin-halasi}

\subsection{Main result} One of the main algebraic results about representations of
finite algebra groups is the following

\begin{thm}[Gutkin-Halasi]\label{t:gutkin-halasi}
Every irreducible representation of a finite algebra group $1+A$
is induced from a one-dimensional representation of a subgroup of
the form $1+B$, where $B\subseteq A$ is a subalgebra.
\end{thm}

It was stated by E.A.~Gutkin in \cite{gutkin}; however, to the
best of our knowledge, the first complete proof was given by
Z.~Halasi in \cite{halasi}. The following statement is an obvious
consequence of the theorem:

\begin{cor}[Isaacs]\label{c:isaacs}
If $G=1+A$ is a finite algebra group over $\bF_q$, the dimension
of every irreducible representation of $G$ is a power of $q$.
\end{cor}

Historically, however, this corollary was proved by I.M.~Isaacs
before Theorem \ref{t:gutkin-halasi}, using nontrivial techniques
from character theory of finite groups; it appears as Theorem A in
\cite{isaacs}. Moreover, Halasi's proof uses Isaacs' result in an
essential way; on the other hand, Halasi's proof is more
``Lie-theoretic'' than ``group-theoretic.'' In this section we
suggest a self-contained proof of Theorem \ref{t:gutkin-halasi}.
It includes Halasi's argument, which we reproduce here for the
reader's convenience, and also because the technique of the proof
will be useful later on. However, our proof avoids Isaacs' Theorem
A by making use of the notion of a strongly Heisenberg
representation (which does not appear in \cite{isaacs} and
\cite{halasi}). The advantage of our approach (apart from the fact
that it yields a shorter proof) is that now the only nontrivial
results from character theory that we use are the well-known
Frobenius reciprocity theorem and Mackey's irreducibility
criterion.

\subsection{Proof of Theorem \ref{t:gutkin-halasi}}
We will use the following result:
\begin{thm}[\cite{halasi}, Theorem 1.3]\label{t:halasi} Let $A$ be a finite dimensional
associative nilpotent algebra over $\bF_q$, and let $G=1+A$. If
$\phi$ is an irreducible character of $1+A^2$ which is
$G$-invariant, then $\phi$ is linear $($i.e., has degree $1${}$)$.
\end{thm}

\mbr

Now, in the situation of Theorem \ref{t:halasi}, let $\rho:G\to
GL(V)$ be an irreducible representation. We may assume that $\dim
V>1$, in which case $\rho\big\lvert_{1+A^2}$ is not irreducible by
Theorem \ref{t:halasi}. Our proof uses induction on $\dim_{\bF_q}
A$ (the case when $\dim_{\bF_q} A=1$ being obvious). Thus we may
further assume that the statement of Theorem
\ref{t:gutkin-halasi}, and hence, a fortiori, the statement of
Isaacs' Theorem A \cite{isaacs}, holds for all finite algebra
groups $1+B$ with $\dim B<\dim A$. Using the transitivity of
induction of group representations, it suffices to show that
$\rho\cong\Ind^G_H\phi$ for some proper algebra subgroup
$H\subsetneq G$ and some representation $\phi$ of $H$.

\mbr

Let $\psi:1+A^2\to GL(W)$ be an irreducible constituent of
$\rho\bigl\lvert_{1+A^2}$, and let $H=1+U\supseteq 1+A^2$ be a
maximal algebra subgroup of $G$ such that $\psi$ can be extended
to $H$. We have, by Frobenius reciprocity,
\[
\bigl[ \Ind_{1+A^2}^H \psi : \Res_H^G\rho \bigr] = \bigl[ \psi :
\rho\bigl\lvert_{1+A^2} \bigr] \geq 1,
\]
whence $\Ind_{1+A^2}^H\psi$ and $\Res_H^G\rho$ share at least one
irreducible summand, call it $\phi$. Since $\psi$ can be extended
to $H$, and since $1+A^2$ is normal in $H$, it follows that $\phi$
is an extension of $\psi$. We also have, by construction,
\[
\bigl[ \rho : \Ind_H^G\phi \bigr] = \bigl[ \Res_H^G\rho : \phi
\bigr] \geq 1.
\]
Hence, if we prove that $\Ind_H^G\phi$ is irreducible, the
argument will be complete.

\mbr

The proof uses Mackey's criterion. Since $H$ is normal in $G$
(because $H\supseteq 1+A^2$ and $(1+A)/(1+A^2)$ is abelian), we
must show that for any $x\in A$ such that $1+x\not\in H$, the
representations $\phi$ and $\phi^x$ are nonisomorphic, where
\[
\phi^x :  H \rar{} GL(W)
\]
is defined by
\[
\phi^x(h) = \phi\bigl( (1+x)h(1+x)^{-1} \bigr).
\]
Fix $x\in A\setminus U$, and set $N_x=1+\bF_q x+U$; this is of
course an algebra subgroup of $G$ such that $[N_x:H]=q$. Pick an
arbitrary irreducible representation $\la$ of $N_x$ such that
$\phi$ is a summand of $\Res_H^{N_x}\la$.

\mbr

Up to this point we have closely followed the proof given in
\cite{halasi}. Here Halasi applies Isaacs' Theorem A, which shows
that both $\dim\phi$ and $\dim\la$ are powers of $q$. On the other
hand,
\[
\dim\phi < \dim\la \leq \dim\bigl( \Ind_H^{N_x}\phi \bigr) =
q\cdot\dim\phi,
\]
where the first inequality is due to the fact that $\phi$ cannot
be extended to $N_x$ by the choice of $H$. Thus the second
inequality must be an equality, whence $\Ind_H^{N_x}\phi\cong\la$
is irreducible. Applying the converse of Mackey's criterion
implies that $\phi^x\not\cong\phi$ and completes the proof.

\mbr

We now explain how the use of Isaacs' Theorem A can be avoided.
Note that if $N_x\neq G$, then by induction we already know that
Isaacs' Theorem A holds for $H$ and for $N_x$, whence there is no
problem. So let us assume that $N_x=G$, or, equivalently,
$\dim(A/U)=1$. We now use an idea that appeared in \cite{isaacs}.

\mbr

Pick an arbitrary subspace $U'\subset A$ such that $U'\neq U$,
$A^2\subseteq U'$, and $\dim(A/U')=1$. Then $H'=1+U'$ is also an
algebra subgroup of $G$; moreover, $HH'=G$. Let $V=U\cap U'$ and
$M=1+V=H\cap H'$; we have $[H':M]=q$. If the representation
$\phi_0:=\phi\big\vert_{M}$ extends to $H'$, then $\phi_0$, and
hence, a fortiori, $\psi$ is invariant under both $H$ and $H'$,
and hence under $G$. Now Theorem \ref{t:halasi} implies that
$\psi$ has dimension $1$. Moreover, we see that
$\rho\bigl\lvert_{1+A^2}$ decomposes into a direct sum of copies
of $\psi$, so $\rho(1+A^2)$ consists of scalars, i.e., $\rho$ is a
strongly Heisenberg representation. But in this case, by
Proposition \ref{p:strongly-Heisenberg} and Lemma
\ref{l:heisenberg}, we already know that Theorem
\ref{t:gutkin-halasi} is valid.

\mbr

The only remaining case is where $\phi_0$ does not extend to $H'$.
In this case, we can apply Isaacs' Theorem A to the groups $M$ and
$H'$, so Halasi's dimension counting argument given above shows
that $\Ind_M^{H'}\phi_0$ is irreducible. Mackey's criterion
implies that the stabilizer of $\phi_0$ in $H'$ is equal to $M$. A
fortiori, if $G^\phi$ denotes the stabilizer of the isomorphism
class of $\phi$ under the action of $G$, then $H\subseteq G^\phi$
and $G^\phi\cap H'=M$. Since $G=HH'$, these two conditions force
$G^\phi=H$, and so $\Ind_H^G\phi$ is irreducible (again by
Mackey's criterion), which finally completes the proof of Theorem
\ref{t:gutkin-halasi}.

\section{Reduction process for finite algebra groups}\label{s:reduction}

\subsection{} In this section we show that any irreducible
representation $\rho:G\to GL(V)$ of an algebra group $G=1+A$ over
a finite field $k\subset\bF$ can be represented as the induced
representation of a strongly Heisenberg representation of an
algebra subgroup in a canonical way. The precise meaning of this
statement is as follows.
\begin{thm}[Reduction process]\label{t:reduction}
There exists a canonical (i.e., independent of any choices)
decomposition of $V$ into a direct sum of subspaces
\begin{equation}\label{e:3.1}
V=\bigoplus_{i\in I} V_i
\end{equation}
satisfying the following properties.
\begin{enumerate}[(1)]
\item For each $g\in G$ and each $i\in I$, there exists $g(i)\in
I$ with $g(V_i)\subseteq V_{g(i)}$.
\item For some, and hence every, $i\in I$, the group
$G_i=\{g\in G \lvert g(i)=i\}$ is an algebra subgroup of $G$.
\item For some, and hence every, $i\in I$, the representation of
$G_i$ in $V_i$ is strongly Heisenberg.
\item If $\rho':G'\to GL(V')$ is an irreducible representation of
another algebra group $G'=1+A'$ over $k$, if $\phi:A\to A'$ is an
algebra isomorphism, and if $\psi:V\to V'$ is a vector space
isomorphism which intertwines $\rho\circ\phi^{-1}$ and $\rho'$,
then $\psi$ takes the canonical decomposition of $V$ onto the
canonical decomposition of $V'$.
\end{enumerate}
\end{thm}
We remark that, in the situation of the theorem, if we choose
$i\in I$, then properties (1) and (2) imply that $\rho$ is
naturally isomorphic to $\Ind_{G_i}^G\rho_i$, where $\rho_i$
denotes the representation of $G_i$ in $V_i$. However, we prefer
not to choose $i$, which yields a more invariant statement.

\subsection{} The proof of the theorem is easy using the
techniques developed in Section \ref{s:gutkin-halasi}. Namely, by
induction on $\dim V$, it suffices to show that if $\rho$ is not
strongly Heisenberg, then there exists a {\em nontrivial}
decomposition \ref{e:3.1} which is independent of any choices and
satisfies properties (1), (2) and (4), but not necessarily (3), of
the theorem.

\mbr

Suppose $\rho$ is not strongly Heisenberg, and let $I$ denote the
set of isomorphism classes of the irreducible constituents of
$\rho\big\vert_{1+A^2}$. From Theorem \ref{t:halasi}, it follows
that $\abs{I}>1$. For each $i\in I$, let $V_i\subset V$ denote the
sum of all irreducible subrepresentations of
$\rho\big\vert_{1+A^2}$ corresponding to $i$. Then it is clear
that we obtain a decomposition of $V$ of the form \eqref{e:3.1}
satisfying properties (1) and (4), and we only need to verify
property (2). It follows from the more general

\begin{lem}\label{l:3.2}
Let $G=1+A$ be an algebra group over $\bF_q$, let $U\subseteq A$
be any subspace such that $A^2\subseteq U$, so that $H=1+U$ is
automatically a normal algebra subgroup of $G$, and let $\psi$ be
any irreducible representation of $H$. Then the stabilizer
$G^\psi$ of (the isomorphism class of) $\psi$ in $G$ is an algebra
subgroup of $G$.
\end{lem}
\begin{proof}
Since $G/H$ is isomorphic to the underlying additive group of
$A/U$, and since $H\subseteq G^\psi$, it suffices to show that if
$1+x\in G^\psi$, then $1+\la x\in G^\psi$ for all $\la\in \bF_q$.
This is clear if $x\in U$, so let $x\in A\setminus U$ be such that
$1+x\in G^\psi$, and let $N_x=1+U+\bF_q\cdot x$ as in the proof of
Theorem \ref{t:gutkin-halasi}. Then $N_x$ is an algebra group, and
since $x\in G^\psi$, Mackey's criterion implies that
$\Ind_H^{N_x}\psi$ is not irreducible.

\mbr

Let $\phi$ be an irreducible summand of $\Ind_H^{N_x}\psi$; then
$\dim(\psi)\leq\dim(\phi)<q\cdot\dim(\psi)$, whence
$\dim\phi=\dim\psi$ by Corollary \ref{c:isaacs}, and therefore
$\psi\cong\phi\big\vert_H$ by Frobenius reciprocity. A fortiori,
we have $N_x\subseteq G^\psi$, which completes the proof of the
lemma.
\end{proof}

\mbr

This also finishes the proof of Theorem \ref{t:reduction}. We
should point out, however, that in practice it will be more
convenient to use not the whole reduction process, but rather its
first step. For instance, it has the advantage of producing
algebra subgroups $G_i\subseteq G$ that contain $1+A^2$, and, in
particular, are normal.

\subsection{} We conclude this section by introducing two
important invariants of irreducible representations of finite
algebra groups. If $G=1+A$ is a finite algebra group over $\bF_q$
and $\rho:G\to GL(V)$ is an irreducible representation, we define
the {\em functional dimension} of $\rho$ to be
$\fdim(\rho)=\log_q(\dim_\bC V)$. This definition is motivated by
the geometric applications; note that $\fdim(\rho)$ is a
nonnegative integer by Corollary \ref{c:isaacs}. On the other
hand, we define $\sh(\rho)$ to be the number of steps in the
reduction process applied to $\rho$, with the understanding that
$\sh(\rho)=0$ if and only if $\rho$ is strongly Heisenberg. We
will see later on (in Theorem \ref{t:base-general}) that both of
these invariants are stable under base change maps.

\section{Base change maps for $1$-dimensional
representations}\label{s:base-1-dim}

\subsection{Notation and terminology}
In this section we begin studying the aspects of character theory
that have to do with extending the base field. All fields we
consider throughout the rest of the paper are assumed to be finite
subfields of $\bF$, without any exceptions. If the ``base field''
is $\bF_q$, and if $k=\bF_{q^m}$ and $k'=\bF_{q^n}$ are such that
$k\subseteq k'$, i.e., $m\lvert n$, then the base change map
$T_m^n$ will also be denoted by $T_{k}^{k'}$. (This is a slight
abuse of notation, since in principle the base change maps may
depend on the choice of the base field $\bF_q$; however, as we
will see, all our constructions depend only on the pair of fields
$k\subseteq k'$.)

\mbr

If $k=\bF_q$, the topological generator $\Fr_q$ of the Galois group
$\Gal(\bF/k)$ will also be denoted by $\Fr_k$. From now on, if $A$
is a finite dimensional associative nilpotent algebra over $k$, we
will say that $G=1+A$ is an {\em algebra group over $k$}, and we
will no longer make an explicit distinction between finite algebra
groups and unipotent algebra groups; we will also write $\dim
G=\dim_k A$. If $k\subseteq k'\subseteq k''$ are field extensions,
we will write $A'=k'\tens_k A$, $A''=k''\tens_k A=k''\tens_{k'}
A'$, and we will write $G'=1+A'$, $G''=1+A''$ for the algebra
groups obtained by extending scalars from $k$ to $k'$ and $k''$,
respectively.

\mbr

Finally, if $\bF_q\subseteq k$ and $G=1+A$ is an algebra group
over $k$, we will say that $A$ or $G$ {\em is defined over
$\bF_q$} if $A=k\tens_{\bF_q} A_0$ for some algebra $A_0$ over
$\bF_q$. We also make the convention that whenever a fixed power
$q$ of the prime $p$ is present, we write $\Fr$ in place of
$\Fr_q$.

\subsection{Norm maps} For any group $G$, we will denote by $G^{ab}$ its
abelianization, i.e., $G/(G,G)$. The main result of this section
is the following

\begin{thm}\label{t:norms}
For all field extensions $k\subseteq k'$, and for all algebra
groups $G=1+A$ over $k$, there exist group homomorphisms
\[
N_{k'/k}^A : G'^{ab}=\bigl( 1+A' \bigr)^{ab} \rar{} \bigl( 1+A
\bigr)^{ab} = G^{ab}
\]
with the following properties:
\begin{enumerate}[(1)]
\item $N^A_{k'/k}$ is functorial with respect to $A$ (in the
obvious sense);
\item if $A$ is defined over $\bF_q\subseteq k$, then $N^A_{k'/k}$
is $\Fr$-equivariant;
\item if $k'\subseteq k''$ is another extension, we have
\[
N^A_{k''/k} = N^A_{k'/k} \circ N^A_{k''/k'};
\]
\item the homomorphism $N^A_{k'/k}$ is surjective.
\end{enumerate}
\end{thm}

The homomorphisms $N_{k'/k}^A$ of the theorem will be called the
{\em norm maps}. Due to the functoriality property of the norm
maps, there is little harm in omitting the superscript $A$ from
the notation, which we will do from now on.

\mbr

Let us now make a few comments about the statement of the theorem,
all of which follow easily from the proof given below. The
construction of $N_{k'/k}$ depends only on $k$ and $k'$, and not
on the field over which $A$ is defined. Property (2) implies in
particular that $N_{k'/k}$ is invariant with respect to $\Fr_k$,
whence the kernel of $N_{k'/k}$ contains all elements of the form
$\Fr_k(g)\cdot g^{-1}$, for $g\in G'$. Property (4) then implies
that $N_{k'/k}$ yields a surjection
\begin{equation}\label{e:norms-h0}
H_0\bigl( \Gal(k'/k), G'^{ab} \bigr) \twoheadrightarrow G^{ab}.
\end{equation}

\subsection{Base change maps} For any group $G$, we will denote by $G^*$ the group of all
homomorphisms $G\to\cst$, i.e., the group of isomorphism classes
of one-dimensional representations of $G$. Of course, we have
$G^*=(G^{ab})^*$. Moreover, since the functor $G\mapsto G^*$ is an
anti-autoequivalence of the category of finite abelian groups, it
is easy to see that the existence of norm maps with properties
(1)--(3) of the theorem is equivalent to the existence of base
change maps for one-dimensional representations of finite algebra
groups, and property (4) is equivalent to the injectivity of the
base change maps, since the latter are precisely the duals of the
homomorphisms \eqref{e:norms-h0}. These remarks will be used
implicitly from now on.

\mbr

As mentioned in the introduction, it is not known to us whether
the base change maps $G^*\to (G'^*)^{\Gal(k'/k)}$ constructed here
are surjective or not. However, since we are dealing with finite
sets, it is clear that our base change maps are surjective if and
only if $\abs{G^{ab}}=\abs{H_0(\Gal(k'/k),G'^{ab})}$ for all
extensions $k\subseteq k'$ as above. Note that since all such
extensions are cyclic, the last equality is equivalent to the
equality
\begin{equation}\label{e:orders}
\bigl\lvert G^{ab} \bigr\rvert = \bigl\lvert(G'^{ab})^{\Gal(k'/k)}
\bigr\rvert.
\end{equation}
On the other hand, if our base change maps are not surjective,
then, for the same reason, no surjective base change maps can
possibly exist.

\mbr

We will see in the next two sections that if the base change maps
we have constructed are surjective for $1$-dimensional
representations, then they are surjective for arbitrary
irreducible representations. Moreover, it is easy to see that our
base change maps are necessarily surjective in some special cases.
On the one hand, if the degree $[k':k]$ is relatively prime to
$p$, then they are surjective because the obstruction to the
injectivity of \eqref{e:norms-h0} lies in a certain cohomology group
which vanishes in this case. On the other hand, if the group
$G=1+A$ has nilpotence class less than $p$, then our base change
maps have to be surjective, because there do exist surjective base
change maps, constructed via the orbit method correspondence. If
$A$ satisfies the stronger condition $A^p=(0)$, then the
exponential map $A\to 1+A$ is a well-defined bijection, and it is
easy to check in this case that our base change maps coincide with
the ones provided by the orbit method. However, if we only assume
that $G$ has nilpotence class less than $p$, then it is not clear
to us how to describe the Lie ring scheme associated to $G$ via
Lazard's construction (\cite{khukhro}, \cite{drinfeld}) in terms
of $A$, so we do not know if our base change maps coincide with
the ones provided by the orbit method or not.

\subsection{Proof of Theorem \ref{t:norms}} The idea of the
construction of the norm maps $N^A_{k'/k}$ was explained to the
author by V.~Drinfeld. Given an algebra group $G=1+A$ over $k$,
let us write $R$ for the $k$-algebra obtained by formally
adjoining a unit to $A$. That is, $R$ is the set of formal
expressions of the form $\la\cdot 1+a$, where $\la\in k$ and $a\in
A$, with the obvious algebra structure. We have an obvious
canonical decomposition $R^\times=k^\times\times G$, which induces
a decomposition $(R^\times)^{ab}=k^\times\times G^{ab}$. If we
write $R'=k'\tens_k R$, then we claim that it suffices to define
norm maps $(R'^\times)^{ab}\to (R^\times)^{ab}$ satisfying the
properties stated in the theorem. Indeed, note that $G'^{ab}$ is a
$p$-group, whereas the order of $k^\times$ is relatively prime to
$p$, so there can be no nontrivial homomorphisms from $G'^{ab}$ to
$k^\times$, and also from $k'^\times$ to $G^{ab}$. 
It follows immediately that norm maps for the groups
$(R^\times)^{ab}$ with properties (1)--(4) induce norm maps for
algebra groups, satisfying the same properties.

\mbr

In the remainder of the proof we will freely use the basic notions
and results in algebraic $K$-theory (more precisely, we only need
the functors $K_0$ and $K_1$). For all the missing explanations we
refer the reader to \cite{srinivas}, Chapter 1, and \cite{weibel},
Chapter III. It is proved in \cite{srinivas}, Example (1.6), that
for any algebra $R$ of the form considered above, the natural map
$GL_1(R)\to K_1(R)$ induces an isomorphism
$(R^\times)^{ab}\rar{\simeq} K_1(R)$. On the other hand, if
$R'=k'\tens_k R$, as above, then $R'$ is a projective (in fact,
free) $R$-module of finite rank, so, as explained in
\cite{weibel}, Chapter III, there exists a natural ``transfer
map'' $K_1(R')\to K_1(R)$. Explicitly, the transfer map is induced
by the homomorphism $GL(R')\to GL(R)$ which comes from the fact
that a matrix of size $n$ over $R'$ naturally gives rise to a
matrix of size $n\cdot[k':k]$ over $R$. It is clear from this
description that the norm maps $(R'^\times)^{ab}\to
(R^\times)^{ab}$ induced by the transfer maps in $K$-theory
satisfy properties (1)--(3) of the theorem, and we only need to
verify property (4).

\mbr

To this end, we use an idea suggested by C.~Weibel (private
communication). If $R=k$, then the statement follows trivially
from the multiplicative version of Hilbert's theorem 90. If $R\neq
k$, we may assume by induction that the result holds for all
algebras of smaller dimension. Let $I=\{a\in A\bigl\lvert
aA=Aa=(0)\}$. This is a nontrivial two-sided ideal of $A$, which
leads to a natural exact sequence in $K$-theory (see, e.g.,
\cite{weibel}, Proposition III.2.3)
\begin{equation}\label{e:exact-K}
K_1(R,I) \rar{} K_1(R) \rar{} K_1(R/I) \rar{} K_0(I)
\end{equation}
Moreover, we clearly have $K_0(I)=(0)$, because $I$ is nilpotent.
Now the transfer maps yield a commutative diagram defining a
morphism from a similar exact sequence
\begin{equation}\label{e:exact-prime}
K_1(R',I') \rar{} K_1(R') \rar{} K_1(R'/I') \rar{} K_0(I')=(0)
\end{equation}
to the sequence \eqref{e:exact-K}, where $I'=k'\tens_k I$. The
transfer map $K_1(R'/I')\to K_1(R/I)$ is surjective by the
induction assumption. On the other hand, the description of
$K_1(R,I)$ due to Vaserstein, explained in \cite{weibel}, Exercise
2.4, implies that the natural map $1+I\to K_1(R,I)$ is an
isomorphism. Thus the transfer map $K_1(R',I')\to K_1(R,I)$
reduces to the natural norm map $1+I'\to 1+I$ for the commutative
unipotent group $1+I$, which is well known to be surjective (cf.
\cite{lusztig}, or \cite{drinfeld}, where the norm maps are called
``trace maps'' in the commutative case). Applying the precise
version of the Five Lemma to the diagram formed by the sequences
\eqref{e:exact-K}, \eqref{e:exact-prime} and the transfer maps,
and taking into account the fact that the map $K_0(I')\to K_0(I)$
is injective for trivial reasons, we see that the transfer map
$K_1(R')\to K_1(R)$ is surjective, completing the proof of Theorem
\ref{t:norms}.

\section{Base change maps for strongly Heisenberg
representations}\label{s:base-s-h}

\subsection{} In this section, using the norm maps for unipotent
algebra groups constructed in Section \ref{s:base-1-dim}, we will
define base change maps for strongly Heisenberg representations,
and prove their injectivity. If $G=1+A$ is a finite algebra group,
we will write $\sheis(G)$ for the set of isomorphism classes of
strongly Heisenberg representations of $G$.

\mbr

From now on we will use the following notation. If $\rho:1+A_1\to
GL(V)$ is a representation of a finite algebra group over $k$ and
$f:A_1\to A_2$ is an isomorphism of algebras over $k$, then
$\rho\circ f^{-1}$ is a representation of $1+A_2$, which will be
denoted by $\rho^f$. In particular, suppose that $1+A_1$ is a
normal algebra subgroup of a bigger algebra group $1+A$. Then the
conjugation by any element $g\in 1+A$ induces an algebra
automorphism of $1+A_1$ whose action on representations of $1+A_1$
will be denoted by $\rho\mapsto\rho^g$. The main result of this
section is

\begin{thm}\label{t:base-s-h}
For all field extensions $k\subseteq k'$, and for all algebra
groups $G=1+A$ over $k$, there exist maps
\[
T_k^{k'} : \sheis(G)\rar{}\sheis(G')^{\Gal(k'/k)}
\]
with the following properties:

 \sbr

\begin{enumerate}[(1)]
\item $T_k^{k'}$ is functorial with respect to isomorphisms of algebras over $k$, i.e.,
if $f:A_1\to A_2$ is a $k$-algebra isomorphism and $f':A'_1\to
A'_2$ is obtained from $f$ by extension of scalars, then
$T_k^{k'}(\rho^f)=T_k^{k'}(\rho)^{f'}$ for all
$\rho\in\sheis(1+A)$; \sbr
\item if $A$ is defined over $\bF_q\subseteq k$, then $T_k^{k'}$
is $\Fr$-equivariant, where $\Fr=\Fr_q$; \sbr
\item if $k'\subseteq k''$ is another extension, we have
\[
T_k^{k''} = T_{k'}^{k''} \circ T_k^{k'};
\]
 \sbr
\item the map $T_k^{k'}$ is injective;
\item we have $\fdim(T_k^{k'}(\rho))=\fdim(\rho)$ for all
$\rho\in\sheis(G)$;
\item the map $T_k^{k'}$ is compatible with restriction of
representations in the following sense. Suppose $K=1+U$ is an
algebra subgroup of $G$ such that $1+A^2\subseteq K$, let
$\rho\in\sheis(G)$, and let $\psi\in\widehat{K}$ be an irreducible
summand of $\rho\big\vert_{K}$. Then $\psi$ is also strongly
Heisenberg, and $T_k^{k'}(\psi)\in\sheis(K')$ is an irreducible
summand of $T_k^{k'}(\rho)\big\vert_{K'}$.
\end{enumerate}
\end{thm}

\mbr

Moreover, we will also prove
\begin{prop}\label{p:surj-sh}
If the base change maps for $1$-dimensional representations defined in Section \ref{s:base-1-dim} are surjective for all algebra subgroups of $G$, then the base change maps for strongly Heisenberg representations of $G$ are also surjective.
\end{prop}

\subsection{} To prove Theorem \ref{t:base-s-h}, we need to understand the structure of strongly Heisenberg representations. To this end, let $G=1+A$ be a finite algebra group over $k$, and let $\phi:1+A^2\to\cst$ be a $G$-invariant homomorphism. Lemma \ref{l:strongly-Heisenberg} implies in particular that the kernel of the commutator pairing $c_\phi$ introduced in that lemma is a $k$-subspace of $A/A^2$, which implies that the group
\[
G_\phi = \bigl\{ g\in G \,\bigl\lvert\, \phi(ghg^{-1}h^{-1})=1\ \forall\,h\in G\bigr\}
\]
is an algebra subgroup of $G$. It is not to be confused with the stabilizer $G^\phi$ of $\phi$ in $G$, which in our situation coincides with $G$. Moreover, $c_\phi$ induces a nondegenerate alternating bilinear form on the finite abelian group $G/G_\phi$, and, as explained in the proof of Proposition \ref{p:strongly-Heisenberg}, there exists a Lagrangian subgroup of $G/G_\phi$ which is also a $k$-subspace, i.e., has the form $H/G_\phi$ for an algebra subgroup $H\subseteq G$. By definition, this implies that $\phi$ is trivial on the commutator $(H,H)$, whence $\phi$ admits a (non-unique) extension to a homomorphism $\psi:H\to\cst$. It is then an easy exercise to check that $\Ind_H^G\psi$ is irreducible, and in fact is a strongly Heisenberg representation of $G$. Moreover, it follows from Proposition \ref{p:strongly-Heisenberg} that every strongly Heisenberg representation of $G$ is obtained in this way; more precisely, we have the following result, which follows immediately from the preceeding discussion and from the last statement of Lemma \ref{l:heisenberg}:
\begin{lem}\label{l:classification-sh}
There is a natural bijection between the set of isomorphism classes of strongly Heisenberg representations of $G$ whose restriction to $1+A^2$ is a multiple of $\phi$ and the set of extensions of $\phi$ to a homomorphism $\chi:G_\phi\to\cst$.
\end{lem}
Note that the definition of $G_\phi$ makes sense whenever we have a finite group $G$, a normal subgroup $K\subseteq G$ such that $G/K$ is abelian, and a $G$-invariant character $\phi:K\to\cst$.
The next result is also completely straightforward.
\begin{lem}\label{l:g-phi}
The following are equivalent for a subgroup $H\subseteq G$ such that $K\subseteq H$:
\begin{enumerate}[(i)]
\item $H\subseteq G_\phi$;
\item $\phi$ can be extended to a $G$-invariant homomorphism $H\to\cst$;
\item $\phi$ extends to a homomorphism $H\to\cst$, and every such extension of $\phi$ is $G$-invariant.
\end{enumerate}
\end{lem}

\subsection{} We begin the proof of Theorem \ref{t:base-s-h}. Let $G$ be a finite algebra group over $k$, let $k\subseteq k'$ be an extension, and let $\rho:G\to GL(V)$ be a strongly Heisenberg representation. By Lemmas \ref{l:classification-sh} and \ref{l:g-phi}, we have the corresponding pair $(\phi,\chi)$, such that the restriction of $\rho$ to $1+A^2$ is a multiple of $\phi$, and $\chi$ is an extension of $\phi$ to $G_\phi$, which is necessarily $G$-invariant. Using the norm maps defined in Section \ref{s:base-1-dim}, we obtain characters $\ti{\phi}=\phi\circ N_{k'/k}:1+A'^2\to\cst$ and $\ti{\chi}=\chi\circ N_{k'/k}:(G_\phi)'\to\cst$. We begin by observing that $\ti{\phi}$ is $G'$-invariant. Indeed, we know by Lemma \ref{l:3.2} that the stabilizer of $\ti{\phi}$ in $G'$ is an algebra subgroup of $G'$. On the other hand, if we let $G$ act on $A$ by conjugation, it follows from the functoriality of the norm maps that the induced action of $G$ on $A'$ leaves $\ti{\phi}$ invariant. Thus, if we view $G$ as a subgroup of $G'$ in the natural way, then $G\subseteq G'^{\ti{\phi}}$. As $G'$ is the smallest algebra subgroup of $G'$ containing $G$, this proves that $\ti{\phi}$ is $G'$-invariant.

\mbr

By the same argument, we see that $\ti{\chi}$ is $G'$-invariant, and the functoriality of the norm maps implies that $\ti{\chi}$ is an extension of $\ti{\phi}$. Thus $(G_\phi)'\subseteq (G')_{\ti{\phi}}$ by Lemma \ref{l:g-phi}. To show that this inclusion is an equality, let us write $G_\phi=1+U$, where $U$ is a subspace, $A^2\subseteq U\subseteq A$. Then $(G_\phi)'=1+U'$. On the other hand, $(G')_{\ti{\phi}}$ is invariant under the action of $\Gal(k'/k)$ on $G'$, and is also an algebra subgroup of $G'$, whence $(G')_{\ti{\phi}}=1+V'$ for some subspace $V$ such that $U\subseteq V\subseteq A$. Suppose, for the sake of contradiction, that $U\neq V$. Then there exists a subspace $W\subseteq V$ such that $U\subsetneq W$ and such that $\phi$ admits an extension $\la:1+W\to\cst$. Hence $\ti{\la}=\la\circ N_{k'/k}:1+W'\to\cst$ is an extension of $\ti{\phi}$. Since $W'\subseteq V'$, it follows from Lemma \ref{l:g-phi} that $\ti{\la}$ is $G'$-invariant, and hence, a fortiori, $G$-invariant. Now the surjectivity and the $G$-equivariance of the norm map $N_{k'/k}:(1+W')^{ab}\to (1+W)^{ab}$ implies that $\la$ is also $G$-invariant. However, by Lemma \ref{l:g-phi}, this contradicts the choice of $W$.

\mbr

The conclusion of the last two paragraphs is that the pair $(\ti{\phi},\ti{\chi})$ satisfies the assumptions of Lemma \ref{l:classification-sh}, and hence gives rise to a strongly Heisenberg representation $\ti{\rho}$ of $G'$, determined uniquely up to isomorphism. It is clear that since $\ti{\phi}$ and $\ti{\chi}$ are $\Gal(k'/k)$-invariant, so is $\rho$, whence we have defined the base change maps $T_k^{k'}$ for strongly Heisenberg representations as stated in Theorem \ref{t:base-s-h}. Property (5) is obvious from the construction, since if $\rho\in\sheis(G)$ corresponds to a pair $(\phi,\chi)$ as in Lemma \ref{l:classification-sh}, then $\fdim(\rho)=\dim_k(G/G_\phi)/2$, and we have $\dim_k(G/G_\phi)=\dim_{k'}(G'/(G')_{\ti{\phi}})$ by the argument given above. Properties (1)--(3) also follow trivially from the construction and the corresponding properties of the norm maps.

\subsection{} Let us prove property (4). Suppose we have two strongly Heisenberg representations, $\rho_1$ and $\rho_2$, of $G$, corresponding to pairs $(\phi_1,\chi_1)$ and $(\phi_2,\chi_2)$, such that $\ti{\rho_1}\cong\ti{\rho_2}$. Using the construction of the base change maps and Lemma \ref{l:classification-sh}, we see immediately that $\ti{\phi_1}=\ti{\phi_2}$, which forces $\phi_1=\phi_2$ by the surjectivity of the norm maps. Then, a fortiori, we have $(G')_{\ti{\phi_1}}=(G')_{\ti{\phi_2}}$, whence $\ti{\rho_1}\cong\ti{\rho_2}$ forces $\ti{\chi_1}=\ti{\chi_2}$ by Lemma \ref{l:classification-sh}, which, in turn, forces $\chi_1=\chi_2$ by the surjectivity of the norm maps.

\mbr

Finally, we prove property (6). Once again, let $(\phi,\chi)$ be the pair of linear characters corresponding to $\rho$. Since $\psi$ is a summand of $\rho\bigl\lvert_K$, it follows that $\phi$ is a summand of $\psi\bigl\lvert_{1+A^2}$. But $\phi$ is invariant under $G$, and hence, a fortiori, under $K$; since $\psi$ is irreducible, we see that $\psi(1+A^2)$ consists only of scalar operators, whence $\psi$ is also strongly Heisenberg, and acts by the character $\phi$ on $1+A^2$. Furthermore, it is clear that $1+A^2\subseteq G_\phi\subseteq K_\phi$. Let $\la$ denote the extension of $\phi$ to $K_\phi$ which corresponds to the representation $\psi$ as in Lemma \ref{l:classification-sh}. Then $\psi\bigl\lvert_{G_\phi}$ is a multiple of $\la\bigl\lvert_{G_\phi}$ and $\rho\bigl\lvert_{G_\phi}$ is a multiple of $\chi$, which implies that $\la$ is an extension of $\chi$. Now the functoriality of the norm maps implies that $\ti{\la}$ is an extension of $\ti{\chi}$, and, reversing the argument above, we see that $\ti{\psi}$ is a summand of $\ti{\rho}\bigl\lvert_{K'}$, which completes the proof of Theorem \ref{t:base-s-h}.

\subsection{} We conclude by proving Proposition \ref{p:surj-sh}. Let $\rho'$ be a $\Gal(k'/k)$-invariant strongly Heisenberg representation of $G'$, and let $(\phi',\chi')$ be the corresponding pair of linear characters. Then both $\phi'$ and $\chi'$ are $\Gal(k'/k)$-invariant. The assumption of the proposition implies that $\phi'=\ti{\phi}$ for some character $\phi:1+A^2\to\cst$. Moreover, the $G$-invariance and surjectivity of the norm maps implies that $\phi$ is $G$-invariant. By the previous argument, $(G')_{\phi'}=(G_\phi)'$. Since $\chi'$ is $\Gal(k'/k)$-invariant, we have $\chi'=\ti{\chi}$ for some character $\chi:G_\phi\to\cst$, by assumption. Moreover, surjectivity of the norm maps implies that $\chi$ is an extension of $\phi$. Now let $\rho$ be the (unique) strongly Heisenberg representation of $G$ defined by the pair $(\phi,\chi)$; we have $\rho'\cong\ti{\rho}$ by construction, which completes the proof.

\section{Base change maps for general irreducible
representations}\label{s:base-general}

\subsection{} We are now ready to define base change maps for
arbitrary irreducible representations of algebra groups. Let
$k\subseteq k'$ be an extension of finite subfields of $\bF$, as
usual; let $G=1+A$ be a finite algebra group over $k$, and let
$\rho:G\to GL(V)$ be an irreducible representation. We define a
representation $\ti{\rho}$ of $G'$, which a priori may not be
irreducible, inductively, as follows. If $\rho$ is strongly
Heisenberg, we set $\ti{\rho}=T_k^{k'}(\rho)$, where $T_k^{k'}$ is
defined in the previous section. If $\rho$ is not strongly
Heisenberg, let
\begin{equation}\label{e:6.1}
V = \bigoplus_{i\in I} V_i
\end{equation}
be the canonical decomposition of $V$ obtained in the first step
of the reduction process of Section \ref{s:reduction}. Pick $i\in
I$, and let $\rho_i$ denote the representation of $G_i$ in $V_i$.
Then $\sh(\rho_i)=\sh(\rho)-1$ (cf. the end of Section
\ref{s:reduction}), so by induction, we may assume that
$\ti{\rho_i}$ has been defined (though it may not be irreducible).
We set
\begin{equation}\label{e:6.2}
\ti{\rho} = \Ind_{G_i'}^{G'}\ti{\rho_i}.
\end{equation}
The main result of this section is the following
\begin{thm}\label{t:base-general}
\begin{enumerate}[(1)]
\item The representation $\ti{\rho}$ defined by \eqref{e:6.2} is
irreducible. \sbr
\item Up to isomorphism, $\ti{\rho}$ is independent of the choice of $i$.
Thus it is legitimate to define
\[
T_k^{k'}(\rho)=\ti{\rho},
\]
 \sbr
\item The map $\rho\mapsto\ti{\rho}$ is functorial with respect to
isomorphisms of algebras over $k$; moreover, if $A$ is defined
over $\bF_q\subseteq k$, then this map commutes with $\Fr_q$. A
fortiori, it commutes with $\Fr_k$, whence, in view of (1) and
(2), we obtain base change maps
\begin{equation}\label{e:base-change}
T_k^{k'} : \Gh \rar{} \bigl(\Gh'\bigr)^{\Gal(k'/k)}, \qquad
\rho\mapsto\ti{\rho},
\end{equation}
which are independent of any choices.
 \sbr
\item Suppose $K\subseteq G$ is an algebra subgroup such that $1+A^2\subseteq K$,
let $\rho\in\Gh$, and let $\psi\in\Kh$ be an irreducible summand
of $\rho\big\vert_{K}$. Then $\ti{\psi}$ is an irreducible summand
of $\ti{\rho}$. \sbr
\item If $k'\subseteq k''$ is another extension, then
\[
T_k^{k''} = T_{k'}^{k''} \circ T_k^{k'}.
\]
\sbr
\item The base change map \eqref{e:base-change} is injective. \sbr
\item We have $\fdim(\ti{\rho})=\fdim(\rho)$ for all $\rho\in\Gh$. \sbr
\item We have $\sh(\ti{\rho})=\sh(\rho)$ for all $\rho\in\Gh$. \sbr
\end{enumerate}
\end{thm}

\mbr

It is clear that this result implies the theorem stated in the
introduction. The reason for stating Theorem \ref{t:base-general}
in this form is that, as will become apparent, the interplay
between the various parts of the theorem is essential for the
inductive proof we present below.

\subsection{} The theorem will be proved simultaneously for all
field extensions $k\subseteq k'\subseteq k''$ by induction on
$\dim G$ and on $\sh(\rho)$. When $\sh(\rho)=0$, it is clear that
all the non-vacuous statements of Theorem \ref{t:base-general}
follow from Theorem \ref{t:base-s-h}, so we may assume that
$\sh(\rho)>0$ and that Theorem \ref{t:base-general} holds for
irreducible representations $\phi$ of arbitrary finite algebra
groups and for arbitrary extensions $k\subseteq k'\subseteq k''$,
as long as $\sh(\phi)<\sh(\rho)$. Similarly, we assume that the theorem holds for all algebra groups of dimension less than $\dim_k A$.

\mbr

We begin by proving parts (1), (2) and (3). Since
$\sh(\rho_i)<\sh(\rho)$, the theorem holds for the group $G_i$ by
assumption, and, in particular, $\ti{\rho_i}$ is irreducible. By
Mackey's criterion, to show that $\ti{\rho}$ is irreducible, we
must check that the stabilizer $G'^{\ti{\rho_i}}$ of (the
isomorphism class of) $\ti{\rho_i}$ in $G'$ equals $G'_i$. By
Lemma \ref{l:3.2}, $G'^{\ti{\rho_i}}$ is an algebra subgroup of
$G'$; moreover, it is clearly $\Gal(k'/k)$-invariant, because
$\ti{\rho_i}$ is $\Gal(k'/k)$-invariant by part (3) of the
theorem. Thus, if $G'^{\ti{\rho_i}}\neq G'_i$, then there exists
an element $g\in G\setminus G_i$ such that $g\in G'^{\ti{\rho_i}}$
(as usual, we view $G$ as a subgroup of $G'$). But conjugation by
$g$ induces an algebra automorphism of the $k$-subalgebra of $A$
corresponding to $G_i$. Since the base change map for $G_i$ is
injective and commutes with the action of $g$ by the induction
assumption, it follows that $g\in G^{\rho_i}$, which contradicts
the irreducibility of $\rho\cong\Ind_{G_i}^G\rho_i$, by Mackey's
criterion. This proves part (1) for the group $G$.

\mbr

Part (2) is easy, since, by construction, for any two elements
$i,j\in I$, there exists an element $g\in G$ which conjugates the pair $(G_i,\rho_i)$ into $(G_j,\rho_j)$. Applying part (3) to $G_i$, we see that $g$ also conjugates $(G'_i,\ti{\rho_i})$ into $(G'_j,\ti{\rho_j})$, whence
\[
\Ind_{G'_i}^{G'} \ti{\rho_i} \cong \Ind_{G'_j}^{G'} \ti{\rho_j},
\]
proving (2). For part (3), we recall that the first step of the reduction process of Section \ref{s:reduction} is functorial with respect to algebra isomorphisms by construction, which implies that if the base change maps for the groups $G_i$ are functorial with respect to algebra isomorphisms, then so are the base change maps for $G$. A similar argument implies that if $A$ is defined over $\bF_q$, then $T_k^{k'}$ commutes with $\Fr=\Fr_q$, which proves (3).

\subsection{} To prove part (4), we first claim that it suffices to consider
the case where $K$ has codimension $1$ in $G$. Indeed, if $K=G$, there is nothing to prove. Otherwise, there exists an algebra subgroup $K_0\subset G$ of codimension $1$ such that $K\subseteq K_0$. Since $\psi$ is a summand of $\rho\bigl\lvert_K$, it follows from Frobenius reciprocity that there exists an irreducible representation $\phi$ of $K_0$ that is a summand of both $\Ind_{K_0}^K\psi$ and $\rho_{K_0}$. If we know part (4) for the codimension $1$ case, then $\ti{\phi}$ is a summand of $\ti{\rho}\bigl\lvert_{K'_0}$. On the other hand, $\psi$ is a summand of $\phi\bigl\lvert_K$ by Frobenius reciprocity, and since $\dim K_0<\dim G$, it follows from the induction assumption that $\ti{\psi}$ is a summand of $\ti{\phi}\bigl\lvert_{K'}$. Combining these two statements yields the desired conclusion.

\mbr

Now we consider the case where $\dim(G/K)=1$. Recall that $\sh(\rho)>0$ by assumption, and consider the decomposition \eqref{e:6.1} obtained in the first step of the reduction process for $\rho$. We fix $i\in I$ and let $\rho_i:G_i\to GL(V_i)$ be the corresponding representation of $G_i$. There are two possibilities. If $K\supseteq G_i$, then $\psi$ is a summand of
\[
\rho\bigl\lvert_K = \Res_K^G \Ind_K^G \bigl( \Ind_{G_i}^K \rho_i \bigr).
\]
Let us write $\phi_i=\Ind_{G_i}^K\rho_i$, which is also irreducible. By Mackey's criterion, $\rho\bigl\lvert_K$ is a direct sum of pairwise nonisomorphic irreducible representations, all of which are conjugate to $\phi_i$. Hence $\psi$ is conjugate to $\phi_i$. But since we are allowed to replace $(G_i,\rho_i)$ with any pair conjugate to it (which is of the form $(G_j,\rho_j)$ for some $j\in I$), and since $G_i$ is normal in $G$ (because it contains $1+A^2$), we may assume, without loss of generality, that $\psi\cong\phi_i$. By the induction assumption, $\ti{\rho_i}$ is a summand of $\ti{\psi}\bigl\lvert_{G'_i}$. By the Frobenius reciprocity,
\[
\bigl[ \ti{\psi} : \ti{\rho}\bigl\lvert_{K'} \bigr] = \bigl[ \Res_{G'_i}^{G'} \Ind_{K'}^{G'} \ti{\psi} : \ti{\rho_i} \bigr],
\]
and since $\Res_{K'}^{G'} \Ind_{K'}^{G'} \ti{\psi}$ certainly contains $\ti{\psi}$ as a summand, we obtain that $[\ti{\psi} : \ti{\rho}\bigl\lvert_{K'}]>0$, as desired.

\mbr

The only remaining case is $K\not\supset G_i$, which forces $G=KG_i$. In this case it is elementary to check that
\[
\Res_K^G \Ind_{G_i}^G \rho_i \cong \Ind_{K\cap G_i}^K \Res^{G_i}_{K\cap G_i} \rho_i,
\]
using only the definition of an induced representations. Thus, Frobenius reciprocity implies that if $\psi$ is a summand of $\rho\bigl\lvert_K$, then
$[\psi\bigl\lvert_{K\cap G_i} : \rho_i\bigl\lvert_{K\cap G_i}]>0$. Let $\phi$ be any irreducible representation of $K\cap G_i$ contained both in $\psi\bigl\lvert_{K\cap G_i}$ and in $\rho_i\bigl\lvert_{K\cap G_i}$. By the induction assumption, it follows that $\ti{\phi}$ is contained both in $\ti{\phi}\bigl\lvert_{K'\cap G_i'}$ and in $\ti{\rho_i}\bigl\lvert_{K'\cap G_i'}$. Reversing the argument above, we find that $\ti{\psi}$ is contained in $\ti{\rho}\bigl\lvert_{K'}$, completing the proof of (4).

\subsection{} Before proceeding with the proof, we need the following
\begin{lem}\label{l:6.4}
In the situation described before the statement of Theorem \ref{t:base-general}, the pair $(G'_i,\ti{\rho_i})$ is among those that appear in the first step of the reduction process applied to the representation $\ti{\rho}$ of $G'$.
\end{lem}
\begin{proof}
By construction, the pair $(G_i,\rho_i)$ is obtained in the following way. Let $\psi_i$ be an irreducible summand of $\rho\bigl\lvert_{1+A^2}$ corresponding to $i\in I$, let $V_i\subseteq V$ be the sum of all subrepresentations of $\rho\bigl\lvert_{1+A^2}$ that are isomorphic to $\psi_i$, let $G_i=G^{\psi_i}$, and let $\rho_i$ denote the representation of $G_i$ in $V_i$. Now by part (4), which we have already proved, $\ti{\psi_i}$ is an irreducible summand of $\ti{\rho}\bigl\lvert_{1+A'^2}$, so we might as well use $\ti{\psi_i}$ in the first step of the reduction process applied to $\ti{\rho}$. On the other hand, we have $G_i\subseteq (G')^{\ti{\psi_i}}$ by part (3) of the theorem, whence $G'_i\subseteq (G')^{\ti{\psi_i}}$ by Lemma \ref{l:3.2}.

\mbr

Conversely, since $(G')^{\ti{\psi_i}}$ is an algebra subgroup of $G'$ and is clearly $\Gal(k'/k)$-stable, we have $G'^{\ti{\psi_i}}=H'$ for some algebra subgroup $H\subseteq G$. Since the base change maps for $1+A^2$ are injective by the induction hypothesis, we deduce that $H\subseteq G^{\psi_i}=G_i$, which implies that $(G')^{\ti{\psi_i}}=G'_i$. On the other hand, we have $\Ind_{G'_i}^{G'}\ti{\rho_i}=\ti{\rho}$ by construction, and $\ti{\rho_i}\bigl\lvert_{1+A'^2}$ contains $\ti{\psi_i}$ as a summand by part (4). Since $\ti{\psi_i}$ is $G'_i$-invariant, also by construction, it follows that $\ti{\rho_i}\bigl\lvert_{1+A'^2}$ is a sum of copies of $\ti{\psi_i}$. Finally, since $\Ind_{G'_i}^{G'}\ti{\rho_i}=\ti{\rho}$ is irreducible, this forces $\ti{\rho_i}$ to be the representation of $G_i$ in the sum of {\em all} irreducible subrepresentations of $\ti{\rho}\bigl\lvert_{1+A'^2}$ that are isomorphic to $\ti{\psi_i}$ and completes the proof of the lemma.
\end{proof}

\subsection{} We proceed with the proof of the theorem. Part (5) follows trivially from the previous lemma and the induction hypothesis. Moreover, Lemma \ref{l:6.4} clearly implies (8) by induction on $\sh(\rho)$, which then implies (7) by induction. It remains to prove part (6). Suppose $\rho,\phi\in\Gh$ are such that $\ti{\rho}\cong\ti{\phi}$. Then $\sh(\rho)=\sh(\phi)$ by part (8), and we will prove that $\rho\cong\phi$ by induction on $\sh(\rho)$, the case $\sh(\rho)=0$ being established in Theorem \ref{t:base-s-h}.

\mbr

Suppose that $\sh(\rho)=\sh(\phi)>0$. Since $\ti{\rho}\cong\ti{\phi}$, and since the first step of the reduction process of Section \ref{s:reduction} is canonically defined, it follows from Lemma \ref{l:6.4} that at least the subgroups $G_i\cong G$ appearing in the first step of the reduction process are the same for $\rho$ and $\phi$. Let $\rho_i$, $\phi_i$ denote the corresponding irreducible representations of $G_i$, so that, in particular, $\rho\cong\Ind_{G_i}^G\rho_i$ and $\phi\cong\Ind_{G_i}^G\phi_i$. By definition,
we now have
\[
\Ind_{G'_i}^{G'}\ti{\rho_i}\cong\ti{\rho}\cong\ti{\phi}\cong \Ind_{G'_i}^{G'}\ti{\phi_i},
\]
which implies, by Frobenius reciprocity, that $\ti{\phi_i}$ is contained in $\ti{\rho}\bigl\lvert_{G'_i}$. However, by Mackey's criterion, $\ti{\rho}\bigl\lvert_{G'_i}$ is a direct sum of pairwise nonisomorphic representations of $G'_i$ that are $G'$-conjugate to $\ti{\rho_i}$. We deduce that there exists a $g\in G'$ which conjugates $\ti{\rho_i}$ into $\ti{\phi_i}$.

\mbr

Recall that $\Fr_k$ denotes the canonical topological generator of $\Gal(\bF/k)$. Since $\ti{\rho_i}$ and $\ti{\phi_i}$ are $\Gal(k'/k)$-invariant, it follows that $\Fr_k(g)$ also conjugates $\ti{\rho_i}$ into $\ti{\phi_i}$. Hence $\Fr_k(g)^{-1}\cdot g\in (G')^{\ti{\rho_i}}=G'_i$. By Lang's theorem, there exists a finite extension $k''\supseteq k'$ and an element $y\in G_i''$ such that $\Fr_k(y)\cdot y^{-1}=\Fr(g^{-1})\cdot g$. We then have $gy\in G$; on the other hand, by construction, $gy$ conjugates $T_k^{k''}(\rho_i)$ into $T_k^{k''}(\phi_i)$. Applying parts (3) and (6) to the group $G_i$, we see that $gy$ conjugates $\rho_i$ into $\phi_i$, which forces $\rho\cong\phi$ and completes the proof.

\subsection{} We now prove the following
\begin{prop}\label{p:surj-general}
If $G=1+A$ is an algebra group over $k$ such that the base change maps for $1$-dimensional representations defined in Section \ref{s:base-1-dim} are surjective for all algebra subgroups of $G$ and all extensions $k'\subseteq k''$ containing $k$, then the base change for arbitrary irreducible representations of $G$ are also surjective.
\end{prop}
\begin{proof}
Let $\phi$ be any $\Gal(k'/k)$-invariant irreducible representation of $G'$. We will prove that $\phi\cong\ti{\rho}$ for some $\rho\in\Gh$ by induction on $\sh(\phi)$, the case $\sh(\phi)=0$ being established in Proposition \ref{p:surj-sh}. Assume that $\sh(\phi)>0$, and let $\psi$ be any irreducible summand of $\phi\bigl\lvert_{1+A'^2}$. Then so is $\Fr_k(\psi)$, whence there exists an element $g\in G'$ such that $\psi^g\cong\Fr_k(\psi)$. Here, for an irreducible representation $\la$ of $1+A'$ and an element $y\in G'$, we denote by $\la^y$ the representation of $1+A'$ obtained by conjugating $\la$ by $y$. By Lang's theorem, there exists a finite extension $k''\supseteq k'$ and an element $x\in G''$ such that $g=\Fr_k(x)\cdot x^{-1}$. A straightforward computation yields
\[
\Fr_k\bigl( T_{k'}^{k''}(\psi)^{x^{-1}} \bigr) \cong T_{k'}^{k''}(\psi)^{x^{-1}}.
\]
Thus $T_{k'}^{k''}(\psi)^{x^{-1}}$ is a $\Fr_k$-invariant summand of $T_{k'}^{k''}(\phi)\bigl\lvert_{1+A''^2}$.

\mbr

Now we may assume, without loss of generality, that $k'=k''$ and that $\psi$ is already $\Fr_k$-invariant. Indeed, if we can prove that $T_{k'}^{k''}(\phi)\cong T_k^{k''}(\rho)$ for some $\rho\in\Gh$, then the injectivity of the base change maps, proved in Theorem \ref{t:base-general}, will imply that $\phi\cong T_k^{k'}(\rho)$. With this assumption, let us apply the first step of the reduction process to $\phi$ and this particular choice of $\psi$. It clearly yields a pair of the form $(G_i',\phi_i)$, where $G_i'=G'^{\psi}$ is obtain by extension of scalars from an algebra subgroup $G_i\subseteq G$. Moreover, since $G_i$ is normal in $G$ (and hence is the same for all $i$), and all possible choices for $\phi_i$ are $G'$-conjugate to each other, the same argument as above using Lang's theorem shows that we may assume that $\phi_i$ is also $\Gal(k'/k)$-invariant, possibly after replacing $k'$ with a finite extension $k''\supseteq k'$.

\mbr

In this case, since $\sh(\phi_i)=\sh(\phi)-1$, we can write $\phi_i=\ti{\rho_i}$ for some $\rho_i\in\Gh_i$ by the induction assumption. The standard argument using Mackey's criterion implies then that $\rho:=\Ind_{G_i}^G\rho_i$ is irreducible. In addition, part (4) of the theorem implies that $\phi_i=\ti{\rho_i}$ is a summand of $\ti{\rho}\bigl\lvert_{G_i'}$, whence $\phi\cong\ti{\rho}$ by Frobenius reciprocity, completing the proof.
\end{proof}

\mbr

This result has the following curious consequence, whose conclusion is independent of the construction of base change maps, but which we were unable to verify directly. Of course, if one can find an example where the conclusion of this corollary does not hold, it would imply that the base change maps for $1$-dimensional representations of algebra groups are not necessarily surjective. The proof is completely straightforward and therefore omitted.
\begin{cor}\label{c:Fr-invariant-restrictions}
Let $G=1+A$ be an algebra group over $k$ satisfying the assumption of Proposition \ref{p:surj-general}, let $H\subseteq G$ be an algebra subgroup such that $1+A^2\subseteq H$, and let $k\subseteq k'$ be an extension. If $\phi\in\Gh'$ is $\Gal(k'/k)$-invariant, then $\phi\bigl\lvert_{H'}$ has a $\Gal(k'/k)$-invariant irreducible summand.
\end{cor}

\end{document}